\documentclass{amsart}
\usepackage{ifthen}
\usepackage{amsmath}
\usepackage{amssymb}    
\usepackage{stmaryrd} 
\usepackage{eufrak}     
\usepackage{mathrsfs} 
\usepackage{bbold}    
\provideboolean{usemathrsfs}
\setboolean{usemathrsfs}{true}
\usepackage{enumerate}
\usepackage{xspace}
\usepackage{boxedminipage}
\usepackage[utf8]{inputenc}
\usepackage{verbatim}
\provideboolean{usebeamer}
\provideboolean{isthesis}
\provideboolean{isamsltex}
\provideboolean{issiamltex}

\usepackage{hyperref}

\newcommand{\Ignore}[1]{}
\newcommand{\freeze}[1]{}
\newcommand{\crossout}[1]{{\textcolor{red!20}{#1}}}
\newcommand{\highlight}[1]{{\textcolor{blue}{#1}}}
\newcommand{\standout}[1]{{\textcolor{magenta}{#1}}}
\provideboolean{shownotes}
\setboolean{shownotes}{true}
\newcounter{margnote}[page]
\newcommand{\margnotemark}{\ensuremath{\text{\standout{\upshape\texttt{>\arabic{margnote}<}}}}}
\newcommand{\margnote}[2][]{
  \ifthenelse{
    \boolean{shownotes}
  }{
    \stepcounter{margnote}
    \margnotemark
    \marginpar{
      \texttt{
        \begin{minipage}{2cm}
          \raggedright\tiny
          \margnotemark{#1}: 
          #2
        \end{minipage}
      }
    }
  }{
  }
}
\newcommand{\mathnote}[2][]{
  \ifthenelse{
    \boolean{shownotes}
  }{
    \stepcounter{margnote}
    \margnotemark
    \text{
      \standout{
        \texttt{
          \tiny
          #2
        }
      }
    }
  }{
  }
}
\provideboolean{showtodo}

\newcommand{\todo}[1]
{\ifthenelse{\boolean{showtodo}}{\margnote[To do.]{#1}}{}}
{\ifthenelse{\boolean{showtodo}}{\end{boxedminipage}}{}}

\provideboolean{showcomments}
\newcommand{\margincomment}[1]{
\ifthenelse{\boolean{showcomments}}{\marginpar{\tiny #1}}{}
}
\provideboolean{showchanges}
\setboolean{showchanges}{false}
\newcommand{\changes}[1]{
  \ifthenelse{\boolean{showchanges}} {{\highlight{#1}}} {#1}
}
\newcommand{\changefromto}[3][replace with]{
  \ifthenelse{\boolean{showchanges}}
  {{\crossout{#2}\margnote{#1}}{\highlight{#3}}}
  {#3\xspace}
}
\newcommand{\ChangePar}[2]{
  \ifthenelse{\boolean{showchanges}}
  {{\par$\mapsfrom$ \textcolor{red!20}{#1}}{\par$\mapsto$ \textcolor{blue}{#2}}}
  {\par #2}
}
\newcommand{\InsertPar}[1]{
  \ifthenelse{\boolean{showchanges}}
  {{\par$\mapsto$ \textcolor{blue}{#1}}}
  {\par #1}
}

\newcommand{\aposteriori}{{aposteriori}\xspace}
\newcommand{\Aposteriori}{{Aposteriori}\xspace}

\usepackage{mathrsfs}
\provideboolean{usemathrsfs}
\setboolean{usemathrsfs}{true}
\newcommand{\mathscript}
	   {\mathscr}

 \newcommand{\qp}[1]{\ensuremath{\!\left({#1}\right)}}

 \newcommand{\powqp}[2]{\ensuremath{\qp{#2}^{\kern -.2em\lower .3ex\hbox{\scriptsize $#1$}}\kern-.3em}}

 \newcommand{\jump}[1]{\ensuremath{\left\llbracket #1\right\rrbracket}}

\newcommand{\registered}%
    {\ensuremath{{}^{\bigcirc\!\;\!\!\!\!\!\!\!\;\text{\sc r}}}}

\newcommand{\diam}{\operatorname{diam}}

\newcommand{\esssup}{\operatorname{ess\,sup}}         

\newcommand{\Eye}[1]{
  \begin{bmatrix}
  \ifthenelse{#1>1}{
    \ifthenelse{#1>2}{
      \ifthenelse{#1>3}{
        1&0&\dotso&0
        \\
        0&1&\dotso&0
        \\
        \vdots&\vdots&\ddots&\vdots
        \\
        0&0&\dotso&1
      }{
        1&0&0
        \\
        0&1&0
        \\
        0&0&1
      }
    }{
      1&0
      \\
      0&1
    }
  }{
    1
  }
  \end{bmatrix}
}
\usepackage{amscd}

\providecommand{\ListParameters}{}
\renewcommand{\ListParameters}
{
	 \setlength{\topsep}{0em}
	 \setlength{\leftmargin}{0em}
         \setlength{\itemsep}{0ex}
	 \setlength{\parsep}{.5ex}
	 \setlength{\itemindent}{\labelsep}
	 \addtolength{\itemindent}{\labelwidth}
}

\newcounter{LetterListItem}
\renewcommand{\theLetterListItem}{(\alph{LetterListItem})}
{
	\begin{list}%
	{\theLetterListItem\ }%
	{\usecounter{LetterListItem}
	 \ListParameters
	}
}%
{\end{list}}
\newcounter{NumberListItem}
\renewcommand{\theNumberListItem}{\arabic{NumberListItem}}
{
	\begin{list}%
	{\theNumberListItem.\ }%
	{\usecounter{NumberListItem}%
	 \ListParameters
	}
}%
{\end{list}}
\newcounter{QuestionListItem}
\renewcommand{\theQuestionListItem}{\textbf{Question \arabic{QuestionListItem}}}
{
	\begin{list}%
	{\theQuestionListItem.\ }%
	{\usecounter{QuestionListItem}%
	 \ListParameters
	}
}%
{\end{list}}
\newcounter{RomanListItem}
\renewcommand{\theRomanListItem}{(\roman{RomanListItem})}
{
	\begin{list}%
	{\theRomanListItem\ }%
	{\usecounter{RomanListItem}
	 \ListParameters
	}
}%
{\end{list}}
\newcounter{StepsItem}
{
	\begin{list}%
	{Step \theStepsItem.\ }%
	{\usecounter{StepsItem}%
	 \ListParameters
	}
}%
{\end{list}}

\usepackage{amsthm} 
\ifthenelse{\boolean{isthesis}}{%
  \setcounter{secnumdepth}{1}%
}{
  \setcounter{secnumdepth}{2} 
}
\providecommand{\ListParameters}{}
\renewcommand{\ListParameters}
{
	 \setlength{\topsep}{0em}
	 \setlength{\leftmargin}{0em}
         \setlength{\itemsep}{0ex}
	 \setlength{\parsep}{.5ex}
	 \setlength{\itemindent}{\labelsep}
	 \addtolength{\itemindent}{\labelwidth}
}
\newtheoremstyle{plain}
  {}
  {}
  {\mdseries\slshape}
  {\parindent}
  {\bfseries}
  {.}
  {.5em}
  {}

\newtheoremstyle{note}
  {}
  {}
  {}
  {\parindent}
  {\bfseries}
  {.}
  {.5em}
  {}

\newtheoremstyle{claim}
  {}
  {}
  {\mdseries\slshape}
  {}
  {\bfseries}
  {}
  {.5em}
  {}

\newtheoremstyle{exercise}
  {}
  {}
  {}
  {}
  {\bfseries}
  {.}
  {1em}
  {}

\newtheoremstyle{break}
  {}
  {}
  {}
  {}
  {\bfseries}
  {.}
  {\newline}
  {}

  \newcommand{\ObsName}{Remark}
  \newcommand{\RemName}{Remark}
  \newcommand{\NotName}{Notation}
  \newcommand{\BFNName}{Big~Fat~Note}
  \newcommand{\DefName}{Definition}
  \newcommand{\ExaName}{Example}
  \newcommand{\TheName}{Theorem}
  \newcommand{\LemName}{Lemma}
  \newcommand{\ProName}{Proposition}
  \newcommand{\CorName}{Corollary}
  \newcommand{\PbmName}{Problem}
  \newcommand{\HypName}{Hypothesis}
  \newcommand{\AlgName}{Algorithm}
  \newcommand{\ExeName}{Exercise}
  \newcommand{\SolName}{Solution}
  \newcommand{\ClaName}{Claim}
  \newcommand{\EsyName}{Essay}
  \newcommand{\Proofname}{Proof}

\ifthenelse{\boolean{isthesis}}{%
  \newcommand{\Thecounter}{The}
}{%
  \newcommand{\Thecounter}{subsection}
}

\newcommand{\pdfformat}[1]{
  \provideboolean{pdfoutput}
  \setboolean{pdfoutput}{#1}
  \ifthenelse{\boolean{pdfoutput}}%
	     {\typeout{using pdf}
               \usepackage{pdfsync}
	       \usepackage[pdftex]{graphicx,xcolor}
               \providecommand{\graphext}{pdf}
	       \renewcommand{\graphext}{pdf}
	       \providecommand{\graphextex}{pdf_t}
	       \renewcommand{\graphextex}{pdf_t}
	       \usepackage{epsfig}
	       \usepackage{tikz}
	       \usepackage{rotating}
	     }
	     {
	       \typeout{using eps}
	       \usepackage[dvips]{graphicx,xcolor}
               \providecommand{\graphext}{eps}
	       \renewcommand{\graphext}{eps}
	       \providecommand{\graphextex}{eps_t}
	       \renewcommand{\graphextex}{eps_t}
	       \usepackage{epsfig}
	       \usepackage{tikz}
	       \usepackage{rotating}
	     }
}
\providecommand{\qed}{\vrule height 5pt depth 0pt width 3pt}
\providecommand{\qqed}{{\raggedright{\ \hfill\qed}}}

\newcounter{passo}

\newenvironment{Proof}[1][{}]%
{\par\noindent{\bf \Proofname\ #1}\setcounter{passo}{0}}%
{\qqed\par}
\newenvironment{Proof*}[1][{}]%
{\subsection{\Proofname\ #1}\setcounter{passo}{0}}
{\qqed\par}
\swapnumbers{
  \theoremstyle{plain}
  \ifthenelse
      {\boolean{isthesis}}
      {\newtheorem{The}{\TheName}[section]}
      {}
 {
   \theoremstyle{plain}
   
   \newtheorem{theorem}[\Thecounter]{\TheName}

   \newtheorem{lemma}[\Thecounter]{\LemName}

   \newtheorem*{The*}{\TheName}
   \newtheorem*{Lem*}{\LemName}
   \newtheorem*{Pro*}{\ProName}
   \newtheorem*{Cor*}{\CorName}
   \newtheorem*{Pbm*}{\PbmName}
   \newtheorem*{Hyp*}{\HypName}
   \newtheorem*{Exe*}{\ExeName}
   \newtheorem*{Txx*}{\ExeName} 
   \newtheorem*{Con*}{Conclusion}
   \newtheorem*{Sum*}{Summary}
 }
 {
   \theoremstyle{claim}

 }
 {
   \theoremstyle{note}
   
   \newtheorem{remark}[\Thecounter]{\RemName}
   
   \newtheorem*{Obs*}{\ObsName}

   \newtheorem{definition}[\Thecounter]{\DefName}
   \newtheorem*{Def*}{\DefName}
   
   \newtheorem*{Exa*}{\ExaName}
 }

 {
   \theoremstyle{break}
   
   \newtheorem*{Alg*}{\AlgName}
 }
}

{
  \theoremstyle{exercise}
  \newtheorem{Exe}[subsection]{\ExeName}

  \newtheorem*{Rdn*}{Reading}
  \newtheorem*{Sol*}{\SolName}
}

{\begin{Exe}\vspace{-24pt}\begin{enumerate}[(a)\ ]}%
{\end{enumerate}\end{Exe}}
\renewcommand{\Aposteriori}{{A posteriori}\xspace}
\renewcommand{\aposteriori}{{a posteriori}\xspace}

\newcommand{\be}{\begin{equation}}
\newcommand{\ee}{\end{equation}}
\newcommand{\bea}{\begin{eqnarray}}
\newcommand{\eea}{\end{eqnarray}}
\newcommand{\beaa}{\begin{eqnarray*}}
\newcommand{\eeaa}{\end{eqnarray*}}

\newcommand{\mbf}[1]{\mbox{\boldmath$\rm{#1}$}}

\newcommand{\el}{ \kappa \in \mathcal{T}}

\renewcommand{\diam}{\operatorname{diam}}
\newcommand{\Q}{\mathcal{Q}}
\newcommand{\dint}{\text{\rm int}}

\newcommand{\mean}[1]{\{ \kern -1.6mm \{#1\} \kern -1.6mm \}}
\newcommand{\ha}{\frac{1}{2}}
\renewcommand{\jump}[1]{[ \kern -.7mm [#1] \kern -.7mm ]}

\newcommand{\su}{\sum_{\kappa\in\mathcal{T}}}\newcommand{\ud}{\mathrm{d}}
\newcommand{\ndg}[1]{| \kern -.25mm \|{#1}| \kern -.25mm \|}
\newcommand{\ltwo}[2]{\|{#1}\|_{#2}}

\title[IPDG a posteriori bounds for quasilinear parabolic problems]{A posteriori error bounds for discontinuous Galerkin
methods for quasilinear parabolic problems}
\author{Emmanuil H. Georgoulis}
\address{Department of Mathematics, University of Leicester, University Road,
Leicester LE1 7RH, United Kingdom}
\email{Emmanuil.Georgoulis@mcs.le.ac.uk}
\author{Omar Lakkis}
\address{Department of Mathematics,
University of Sussex, Falmer, East Sussex BN1 9RF, United Kingdom}
\email{O.Lakkis@sussex.ac.uk}
\begin{document}
\begin{abstract}
We derive \aposteriori error bounds for a quasilinear parabolic problem, which is approximated by the
$hp$-version interior penalty discontinuous Galerkin method (IPDG). The error is measured in the energy norm.
The theory is developed for the
semidiscrete case for simplicity, allowing to focus on the
challenges of a posteriori error control of IPDG
space-discretizations of strictly monotone quasilinear parabolic
problems. The a posteriori bounds are derived using the elliptic
reconstruction framework, utilizing available a posteriori error
bounds for the corresponding steady-state elliptic problem.
\end{abstract}
\maketitle

\section{Introduction}
Discontinuous Galerkin (DG) methods \cite{baker,wheeler:78,arnold}, have enjoyed
substantial development in recent years. For parabolic problems DG methods are interesting due to their good local
conservation properties as well as due to their block-diagonal mass matrices.

This work is concerned with the derivation of a posteriori error
bounds for the space-discrete interior penalty discontinuous
Galerkin method (IPDG) for quasilinear parabolic problems with strictly
monotone non-linearities of Lipschitz growth.

A posteriori error bounds for $h$-version DG
methods are derived in \cite{kar-pascal,kar-pascal2,bhl,hsw} and for
DG-in-space parabolic problems in
\cite{sun-wheeler:05,ern-proft:05,yangjiming-chenyaping:06,
georgoulis-lakkis:08,ern-vohralik:09}. The contribution of this work is twofold:
\begin{itemize}
\item the derivation of \aposteriori energy-norm error
bounds for IPDG methods for quasilinear parabolic
problems, and
\item the resulting \aposteriori bounds are are explicit with respect to the local elemental polynomial degree.
\end{itemize}

A key tool in our \aposteriori error analysis is the \emph{elliptic
reconstruction technique} \cite{makridakis-nochetto:03,lakkis-makridakis:06,georgoulis-lakkis:08}.
Roughly speaking, in the elliptic reconstruction framework the
error is split into a \emph{parabolic} and an \emph{elliptic} part, respectively.
In the interest of being explicit with respect
to the dependence of the \aposteriori error bounds in the elemental polynomial degree $p$, we restrict the
presentation to quadrilateral elements of tensor-product type (cf. Remark \ref{remark_quad}).


\section{Model problem and the IPDG method}\label{mod}

Let $\Omega$ be a bounded open (curvilinear) polygonal domain with
Lipschitz boundary $\partial\Omega$ in $\mathbb{R}^d$, $d=2,3$.
For $\omega\subset\Omega$, we consider the standard spaces $L^2(\omega)$ (whose norm is
denoted by $\ltwo{\cdot}{\omega}$ for brevity), $H^1(\omega)$ and $H^1_0(\omega)$,
 whose norm will be denoted by $\ltwo{\cdot}{1}$, along with its dual
$H^{-1}(\Omega)$, with norm $\ltwo{\cdot}{-1}$.
For brevity, the standard inner product on $L^2(\Omega)$ will be
denoted by $\langle \cdot,\cdot\rangle$ and the corresponding norm by $\ltwo{\cdot}{}$.
We also define the spaces $L^2(0,T,X)$, $X\in \{L^2(\omega),H^{\pm 1}(\omega)\}$
and $L^{\infty}(0,T,L^2(\Omega))$,
consisting of all measurable functions $v: [0,T]\to X$, for which
$\|v\|_{L^2(0,T;X)}:=\big(\int_0^T \|v(t)\|_{X}^2\big)^{1/2}<+\infty$ and
$\|v\|_{L^{\infty}(0,T;L^2(\Omega))}:=\esssup_{t\in[0,T]} \|v(t)\|$. (The
differentials in the integrals with respect to $t$
are suppressed for brevity throughout this work.)

We identify function $v\in [0,T]\times\Omega\to \mathbb{R}$ with $v:t\to X$ and we
denote $v(t)$, $t\in[0,T]$, for $v\in [0,T]\times\Omega\to \mathbb{R}$.

For $t\in(0,T]$, we consider the problem of finding a function $u$ satisfying
\begin{equation}\label{pde}
 u_t(t,x)-\nabla \cdot (a(t,x, |\nabla u(t,x)|)\nabla u(t,x))=f(t,x)\quad\text{in } (0,T]\times\Omega,
\end{equation}
where $f\in L^{\infty}(0,T; L^2(\Omega))$ and
$a$ scalar uniformly continuous function, subject to initial
condition
 $ u(0,x)=u_0(x)$ on $\{0\}\times\Omega$,
for $u_0\in L^2(\Omega)$, and homogeneous Dirichlet boundary
conditions on $[0,T]\times\partial\Omega$.

We assume that the non-linearity
$a$ in equation (\ref{pde}) is of strongly monotone type with Lipschitz growth so that
there exist positive constants $\underline{a}$ and $\overline{a}$ such that the following inequalities hold:
\bea
|a(t,x,|y|)y-a(t,x,|z|)z|&\le& \overline{a}|y-z|\label{growth_above}\\
\big(a(t,x,|y|)y-a(t,x,|z|)z\big)\cdot(y-z)&\ge& \underline{a}|y-z|^2\label{growth_below},
\eea
for all vectors $y,z\in \mathbb{R}^d$, and all $(t,x)\in [0,T]\times\bar{\Omega}$.

Let $\mathcal{T}$ be a shape-regular subdivision of $\Omega$ into disjoint closed quadrilateral
elements $\el$. We assume that $\el$ are
constructed via $\mathrm{C}^{\infty}$-diffeomorphisms with non-singular Jacobian
$F_{\kappa}:(-1,1)^d\to\kappa$, so as to ensure $\bar{\Omega}=\cup_{\el}\bar{\kappa}$.

For $p\in\mathbb{N}$, $\Q_{p}(\hat{\kappa})$ is
the set of all tensor-product polynomials on $(-1,1)^d$ of
degree $p$ in each variable and let
\begin{equation}\label{eq:FEM-spc}
S^p:=\{v\in L^2(\Omega):v|_{F_{\kappa}}
    \in \Q_{p}(\hat{\kappa}),\,\el\},
\end{equation}
be the (discontinuous) \emph{finite element space}.
Let $\Gamma$ be the union of all $(d-1)$-dimensional element faces $e$
associated with the subdivision $\mathcal{T}$ (including the
boundary). Let also $\Gamma_{\dint}:=\Gamma\backslash\partial\Omega$, so that
$\Gamma=\partial\Omega\cup\Gamma_{\dint}$.

Let $\kappa^+$, $\kappa^-$ be two (generic) elements sharing a face
$e:=\kappa^+\cap\kappa^-\subset\Gamma_{\dint}$ with
respective outward normal unit vectors $\mbf{n}^+$ and $\mbf{n}^-$ on $e$.
For $q:\Omega\to\mathbb{R}$ and
$\mbf{\phi}:\Omega\to\mathbb{R}^d$, let
$q^{\pm}:=q|_{e\cap\partial\kappa^{\pm}}$ and
$\mbf{\phi}^{\pm}:=\mbf{\phi}|_{e\cap\partial\kappa^{\pm}}$,
and set
\[
\begin{aligned}
& \mean{q}|_e:=\ha(q^+ + q^-),\ & \mean{\mbf{\phi}}|_e:=\ha(\phi^+ + \phi^-),\\
& \jump{q}|_e:=q^+\mbf{n}^++q^-\mbf{n}^-,\ &
\jump{\mbf{\phi}}|_e:=\phi^+\cdot \mbf{n}^++\mbf{\phi}^-\cdot \mbf{n}^-;
\end{aligned}
\] if $e\subset \partial\kappa\cap\partial\Omega$, we set
$\mean{\mbf{\phi}}|_e:=\mbf{\phi}^+$  and $\jump{q}|_e:=q^+\mbf{n}^+$.
Finally, we introduce the \emph{meshsize}
$h:\Omega\to \mathbb{R}$, defined by $h(x)= \diam{\kappa}$,
if $x\in \kappa\backslash\partial\kappa$ and $h(x)= \mean{ h}$, if $x\in \Gamma$.

Consider the IPDG semi-linear form
$B(\cdot,\cdot):S^p\times S^p\to \mathbb{R}$, introduced in \cite{hrs} for the solution of the
corresponding steady-state problem, defined by
\begin{equation}\label{dgsemilinear}
\begin{aligned}
 B(w,v):=&\su{}\int_{\kappa}\alpha(w)\cdot\nabla v\,\ud x
       +\int_{\Gamma}\big(\theta\mean{a(t,x,h^{-1} |\jump{w}|)\nabla v}\cdot\jump{w}\\
       &-\mean{\alpha(w)}\cdot\jump{v}+\sigma\jump{w}\cdot\jump{v}\big)\,\ud s,
\end{aligned}
\end{equation}
where $\alpha(w):=a(t,\cdot,|\nabla w|)\nabla w$, $w\in H^1(\Omega)+S^p$,
for $\theta\in\{-1,0,1\}$, with the function $\sigma:\Gamma\to \mathbb{R}_+$ defined piecewise by
$\sigma|_e:=C_{\sigma} p^2/(h|_e)$,
for some sufficient large constant $C_{\sigma}>0$.
The corresponding energy norm $\ndg{\cdot}$ is defined
$  \ndg{w}:=\left(\su\ltwo{\nabla w}{\kappa}^2
    +\int_{\Gamma}\sigma\jump{w}^2\ud s\right)^{1/2}$, for $w\in H^1(\Omega)+S^p$.
The (spatially semidiscrete) \emph{interior penalty discontinuous Galerkin method} (IPDG)
for the initial/boundary value model problem reads:
 \be\label{dg_semi}
\text{find}\ U:(0,T]\to S^p\text{ such that  }\\
\langle  U_t,V\rangle+B(U,V)=\langle f,V\rangle \ \forall t\in(0,T], V\in S^p.
\ee

\section{\Aposteriori error bounds}\label{sec:semi-discrete}
For $w\in H^1(\Omega)+S^p$, and $T>0$, we define the norm
$
\ndg{w}_{L^2(0,T;H^1(\Omega))}:=\left(\int_0^T\ndg{w}^2\right)^{1/2},
$
$t>0$. We shall derive \aposteriori bounds for the error
$\ndg{u-U}_{L^2(0,T;H^1(\Omega))}$.

\begin{definition}[elliptic reconstruction]\label{elliptic_recon}
Let $U$ be the (semi-discrete) solution to the problem (\ref{dg_semi}) and fix $t\in[0,T]$.
We define the \emph{elliptic reconstruction} $w\equiv w(t)\in H^1_0(\Omega)$ of $U$ to be
the solution to the elliptic problem
\be\label{ell_rec}
\langle \alpha(w), \nabla v\rangle
= \langle g, v\rangle\quad \forall v\in H^1_0(\Omega),
\ee
where $g\equiv g(t)$ is given by
$g:=-AU+f-\Pi f$, with $\Pi:L^2(\Omega)\to S^p$ is the orthogonal $L^2$-projection operator onto $S^p$
and $A\equiv A(t):S^p\to S^p$ is the discrete operator defined by
\be\label{ell_rec2}
\text{for}\ Z\in S^p,\quad \langle -A Z, V \rangle = B(Z,V)\quad \forall V\in S^p.
\ee
\end{definition}
The construction of $w$ and that of $AZ$ are both well defined in view of the elliptic problem's
unique solvability and the Riesz representation, respectively.
\begin{remark}\label{remark_recon}
The key property of the construction in Definition \ref{elliptic_recon} is that $U$
is the IPDG solution of an elliptic problem with analytical solution $w$. Namely, for each fixed
$t\in[0,T]$ it satisfies
\be\label{dg_steady}
\text{find}\ U\in S^p\text{ such that  } B(U,V)=\langle g,V\rangle
\quad \forall V\in S^p.
\ee
\end{remark}

We can now decompose the error as follows:
\be\label{splitting_semi}
U-u=\rho-\epsilon ,\ \text{with}\ \rho:=w-u,\ \text{and}\ \epsilon:=w-U,
\ee
where $w\equiv w(t)$ denotes the elliptic reconstruction of $U\equiv U(t)$, $t\in[0,T]$.

\begin{lemma}[differential error relation]\label{energy_arg} Let $u$, $w$, $U$, $e$, $\rho$, $\epsilon$ as above. Then,
for all $v\in H^1_0(\Omega)$, we have
\be\label{en_ar}
\langle e_t,v\rangle +\langle \alpha(w)-\alpha(u),\nabla v\rangle =0.
\ee
\end{lemma}
\begin{Proof} We have
\begin{equation}\label{rho-term}
\begin{aligned}
&\langle e_t,v\rangle +\langle \alpha(w)-\alpha(u),\nabla v\rangle
 = \langle U_t,v\rangle +\langle \alpha(w),\nabla v\rangle - \langle f, v\rangle \\
& = \langle U_t,v\rangle +\langle g, v\rangle - \langle f, v\rangle
 = \langle U_t,v\rangle +\langle -AU, v\rangle - \langle \Pi f, v\rangle \\
& = \langle U_t,\Pi v\rangle +\langle -AU, \Pi v\rangle - \langle  f, \Pi v\rangle =0,
\end{aligned}
\end{equation}
using (\ref{pde}), (\ref{ell_rec}) and the properties of the $L^2$-projection, respectively.
\end{Proof}

We consider further the decomposition of $U$ into \emph{conforming} and \emph{non-conforming} (discontinuous) parts
$U = U^c+U^d$,
where $U^c\in S^p\cap H^1_0(\Omega)$ and $U^d:=U-U^c\in S^p$.
Note that there are many ways of performing this decomposition (e.g., by projecting $U$ onto the conforming space)
whereof the specific nature remains at our disposal until further.

We also use the shorthand notation
$e^c := U^c-u$ and $\epsilon^c:=w-U^c$; note that $e^c=\rho-\epsilon^c$,
$e=e^c+U^d$ and that $e^c\in H^1_0(\Omega)$.

\begin{theorem}[abstract \aposteriori energy-error estimate]\label{apost_semi_thm}
With $u$, $U$, $U^d$, $e$, and $\epsilon$ as defined above, the following error estimate is satisfied:
\begin{equation}\label{apost_semi}
\begin{aligned}
\ndg{e}_{L^2(0,T,H^1(\Omega))} \le &
C_1\ndg{\epsilon}_{L^2(0,T,H^1(\Omega))}+
 \underline{a}^{-\ha}\big(\ltwo{u_0-U(0)}{}+\ltwo{U^d(0)}{}\big) \\
&
+ C_1\ndg{U^d}_{L^2(0,T,H^1(\Omega))}
+ C_2\ltwo{U^d_t}{L^2(0,T,H^{-1}(\Omega))},
\end{aligned}
\end{equation}
with $C_1:=1+\sqrt{2}\overline{a}\underline{a}^{-1}$ and $C_2:=\sqrt{2}\underline{a}^{-1}$.
\end{theorem}

\begin{Proof} Set $v=e^c$ in (\ref{en_ar}), to deduce
\begin{equation}\label{energy_one}
\langle e^c_t,e^c\rangle +\langle \alpha(U^c)-\alpha(u),\nabla e^c\rangle
=-\langle U^d_t,e^c\rangle +\langle \alpha(U^c)-\alpha(w),\nabla e^c\rangle.
\end{equation}
Conditions (\ref{growth_below}) and (\ref{growth_above}) imply, respectively,
\[
\langle \alpha(U^c)-\alpha(u),\nabla e^c\rangle \ge \underline{a} \ltwo{\nabla e^c}{}^2
\text{, and }\langle \alpha(U^c)-\alpha(w),\nabla e^c\rangle \le \overline{a}
\ltwo{\nabla \epsilon^c}{}\ltwo{\nabla e^c}{},
\]
and the duality pairing $(H^{-1},H^1_0)$ gives
$|\langle U^d_t,e^c\rangle| \le \ltwo{U^d_t}{-1}\ltwo{\nabla e^c}{}.
$
Using the last 3 relations on (\ref{energy_one}), we deduce
\begin{equation}\label{energy_two}
\langle e^c_t,e^c\rangle +\underline{a} \ltwo{\nabla e^c}{}^2
\le
\left(\ltwo{U^d_t}{-1}+\overline{a}\ltwo{\nabla \epsilon^c}{}\right)
\ltwo{\nabla e^c}{},
\end{equation}
which, in turn, implies
\begin{equation}\label{energy_three}
\langle e^c_t,e^c\rangle +\frac{\underline{a}}{2} \ltwo{\nabla e^c}{}^2
\le
\frac{1}{2 \underline{a}}\big(\ltwo{U^d_t}{-1}+\overline{a}\ltwo{\nabla \epsilon^c}{}\big)^2.
\end{equation}
Integrating (\ref{energy_three}) with respect to $t$ between $0$ and $T$, yields
\[
\ltwo{e^c(t)}{}^2 + \underline{a} \int_0^T  \ltwo{\nabla e^c}{}^2
\le \ltwo{e^c(0)}{}^2
+ \frac{1}{\underline{a}}\!\int_0^T\! \left(\ltwo{U^d_t}{-1}+\overline{a}\ltwo{\nabla \epsilon^c}{}\right)^2
,
\]
or
\begin{equation}\label{energy_four}
\begin{aligned}
 \Big(\int_0^T\!  \ltwo{\nabla e^c}{}^2\Big)^{\ha}
\le &  \underline{a}^{-\ha}\ltwo{e^c(0)}{}
+\underline{a}^{-1} \Big(\int_0^T \left(\ltwo{U^d_t}{-1}+\overline{a}\ltwo{\nabla \epsilon^c}{}\right)^2
\Big)^{\ha}\\
\le &  \underline{a}^{-\ha}\ltwo{e^c(0)}{}
+ C_2\ltwo{U^d_t}{L^2(0,T,H^{-1}(\Omega))}\\
&+(C_1-1)\ndg{\epsilon^c}_{L^2(0,T,H^1(\Omega))}
\end{aligned}
\end{equation}
noting that $\ndg{\epsilon^c}=\ltwo{\nabla \epsilon^c}{}$.
Using the bounds $\ndg{\epsilon^c}{}\le \ndg{ \epsilon}{}+\ndg{U^d}{}$,
 $\ltwo{e^c(0)}{} \le \ltwo{e(0)}{}+\ltwo{U^d(0)}{}$  on (\ref{energy_four}) and
the resulting bound on the triangle inequality
\[
\ndg{e}_{L^2(0,T,H^1(\Omega))} \le \ndg{\epsilon^c}_{L^2(0,T,H^1(\Omega))}
+\ndg{U^d}_{L^2(0,T,H^1(\Omega))},
\]
yields the result.
\end{Proof}

For the above result to yield a formally \aposteriori bound, we need to
estimate $\ndg{\epsilon}_{L^2(0,T,H^1(\Omega))}$ further. In particular, in view of Remark \ref{remark_recon},
we require an \aposteriori error bound for the IPDG method for the corresponding elliptic quasilinear
problem (\ref{dg_steady}). Such a result is available in \cite{hsuw}, an instance of which
 and is presented next.

\begin{theorem}[\cite{hsuw}]\label{hp_apost}
Let $w\in H^1_0(\Omega)$ be the elliptic reconstruction defined in
(\ref{ell_rec}) and let $W\in S^p$ be the solution of (\ref{dg_steady}).
Then, for $C_{\sigma}>1$ sufficiently large the bound
\be\label{gen_bound}
\ndg{w-W}^2\le \mathcal{E}(W, g, S^p):=C_{\rm est}\su \Big(\eta_{\kappa}^2+\mathcal{O}(g,W)\Big),
\ee holds, with
\[
\eta_{\kappa}^2=\frac{h_{\kappa}^2}{p^{2}}\ltwo{\tilde{\Pi}\big(g+\nabla\cdot\alpha(W)\big)}{\kappa}^2
+\frac{h_{\kappa}}{p}\ltwo{\tilde{\Pi}_{\Gamma}\jump{\alpha(W)}}{\partial\kappa\backslash\partial\Omega}^2
+C_{\sigma}^2\frac{p^3}{h_{\kappa}} \ltwo{\jump{W}}{\partial\kappa}^2,
\]
and
\[
\mathcal{O}(g,w_{\rm DG})=\su \Big(\frac{h_{\kappa}^2}{p^{2}}
\ltwo{\big(\mathbb{I}-\tilde{\Pi}\big)\big(g+\alpha(W)\big)}{\kappa}^2
+\frac{h_{\kappa}}{p}\ltwo{\big(\mathbb{I}-\tilde{\Pi}_{\Gamma})
\jump{\alpha(W)}}{\partial\kappa\backslash\partial\Omega}^2\Big),
\]
where $\mathbb{I}$ denotes a generic identity operator, $\tilde{\Pi}$ denotes the $L^2$-projection operator onto
$S^{p-1}$, $\tilde{\Pi}_{\Gamma}$ is defined piecewise by $\tilde{\Pi}_{\Gamma}v|_e:= \pi_e^{p-1} v$, for all
elemental faces $e\subset \Gamma$, $v\in L^2(\Omega)$,
where $\pi_e^{p-1}:L^2(\Omega)\to \mathcal{P}_{p-1}(e)$ denotes the $L^2$-projection operator of the trace on the face $e$ of a function in
$S^{p-1}$ (with $\mathcal{P}_{p-1}(e)$, for $e\subset\bar{\kappa}$ the space of mapped univariate polynomials
of degree at most $p-1$ on $e$), and $C_{\rm est}>0$ is
independent of $C_{\sigma}$, $\theta$, $h$ and $p$.
\end{theorem}

Also, it is possible to further estimate the terms involving $U^d$, to avoid computing $U^d$ explicitly.
This is done (with, crucially, \emph{explicit} dependence on $p$)
using the following result based on \cite[Lemma 3.2]{burman-ern}.

\begin{lemma}\label{hp_oswald}
Suppose $\mathcal{T}$ does \emph{not} contain any hanging nodes.
Then, for any $v\in S^p$ and any multi-index $\gamma$, with $|\gamma|=0,1$, there exists
a function $v^c\in S^p\cap H^1_0(\Omega)$ such that
\be\label{magic}
\su \ltwo{D^{\gamma}(v-v^c)}{\kappa}^2\le C_3
\ltwo{\big(\frac{h}{p^2}\big)^{\ha-|\gamma|}\jump{v}}{\Gamma}^2,
\ee
with $C_3>0$ depending on the maximal angle of $\mathcal{T}$ only.
\end{lemma}

\begin{Proof} \cite[Lemma 3.2]{burman-ern} implies that for every $\kappa\in\mathcal{T}$
there exists an Oswald-type operator $I_{\rm Os}:S^p\to S^p\cap H^1_0(\Omega)$, such that
\be
\ltwo{v-I_{\rm Os}v}{\kappa}^2\le C\sum_{e\subset \mathcal{F}(\kappa)} \frac{h_{\kappa}}{p^2}\ltwo{\jump{v}}{e}^2,
\ee
for all $v\in S^p$, with $\mathcal{F}(\kappa):=\{e\in\Gamma:e\cap\bar{\kappa}\neq \emptyset\}$.
Summing over all the elements $\kappa\in\mathcal{T}$, and observing
that the maximal angle and the lack of hanging nodes gives an upper bound on
the cardinality of $\mathcal{F}(\kappa)$ for all $\kappa\in\mathcal{T}$, we deduce that
\be
\su\ltwo{v-I_{\rm Os}v}{\kappa}^2\le C\sum_{e\subset \Gamma} \frac{h_{\kappa}}{p^2}\ltwo{\jump{v}}{e}^2,
\ee
which shows (\ref{magic}) for $|\gamma|=0$.
To show (\ref{magic}) for $|\gamma|=1$, we observe that $(v-I_{\rm Os}v)\in S^p$; thus, the standard inverse
estimate yields:
\be
\su\ltwo{\nabla(v-I_{\rm Os}v)}{\kappa}^2\le C\su\frac{p^4}{h_{\kappa}^2}\ltwo{v-I_{\rm Os}v}{\kappa}^2
\le C\sum_{e\subset \Gamma} \frac{p^2}{h_{\kappa}}\ltwo{\jump{v}}{e}^2,
\ee
using the shape regularity of $\mathcal{T}$. Setting $v^c= I_{\rm Os}v$, the result follows.
\end{Proof}

\begin{remark}\label{remark_quad}
The assumptions of Lemma \ref{hp_oswald} pose the following restrictions on the finite element space $S^p$:
the use of quadrilateral elements (as the tensor-product nature of the local elemental bases is of crucial importance here),
the exclusion of hanging nodes and the uniformity of the polynomial degree.
If explicit knowledge of the polynomial degree $p$ in the \aposteriori bounds presented in this work
is not required, then these restrictions are not needed in view of \cite[Lemma 4.1]{kar-pascal}, i.e.,
triangular elements containing hanging nodes can be employed.
\end{remark}
Combining the results of Theorems \ref{apost_semi_thm} and \ref{hp_apost},
together with the approximation properties described in Lemma
\ref{hp_oswald}, we obtain an \aposteriori error bound in the energy norm for the semi-discrete problem (\ref{dg_semi}).

\begin{theorem}[energy-norm \aposteriori bound]\label{apost_semi_thm_two}
With the notation of Theorem \ref{apost_semi_thm} and the assumptions of Lemma \ref{hp_oswald},
the following error bound holds:
\begin{equation}\label{apost_semi_two}
\begin{aligned}
\ndg{e}_{L^2(0,T,H^1(\Omega))} \le &
C_1\int_0^T \mathcal{E}^2(U, g, S^p)+
 \underline{a}^{-\ha}\ltwo{u_0-U(0)}{}\\
& +\underline{a}^{-\ha}C_3\ltwo{\big(\frac{h}{p^2}\big)^{\ha}\jump{U(0)}}{\Gamma}^2
+ C_4\ltwo{\sqrt{\sigma}\jump{U}}{L^2(0,T;L^2(\Gamma))}\\
&
+ C_5\ltwo{\big(\frac{h}{p^2}\big)^{\ha}\jump{U_t}}{L^2(0,T;L^2(\Gamma))},
\end{aligned}
\end{equation}
with $C_4:=C_1\sqrt{C_3/C_{\sigma}}$, $C_5:=C_2C_{\rm PF}$ and $C_{\rm PF}>0$ (the Poincar\'e--Friedrichs constant),
such that $\ltwo{v}{-1}\le C_{\rm PF}\ltwo{v}{}$, for all $v\in L^2(\Omega)$.
\end{theorem}
\begin{Proof}
Combining the results from Theorems \ref{apost_semi_thm} and \ref{hp_apost}, together with the approximation
properties described in Lemma \ref{hp_oswald}, the result follows.
\end{Proof}
\def\cprime{$'$} \def\cprime{$'$} \def\cprime{$'$} \def\cprime{$'$}
\providecommand{\bysame}{\leavevmode\hbox to3em{\hrulefill}\thinspace}
\providecommand{\MR}{\relax\ifhmode\unskip\space\fi MR }
\providecommand{\MRhref}[2]{%
  \href{http://www.ams.org/mathscinet-getitem?mr=#1}{#2}
}
\providecommand{\href}[2]{#2}

\end{document}